\documentclass [12pt]{article}
\usepackage{amsmath,amsfonts,amssymb}
\usepackage{mathrsfs}

\usepackage[pdftex]{graphicx}

\numberwithin{equation}{section}

\newcommand\qbinom[2]{\genfrac{[}{]}{0pt}{}{#1}{#2}}

\newcommand{\dspace}{\baselineskip 21pt}
\begin{document}

\title{A Common Generalization to Theorems
\\on Set Systems with $\mathcal{L}$-intersections}

\author{
\\Jiuqiang Liu$^{a, b}$, Shenggui Zhang$^{c}$, and Jimeng Xiao$^{c}$
\\
\\$^{a}$ School of Management Engineering
\\Xi'an University of Finance and Economics
\\Xi'an Shaanxi, 710100, PRC
\\
\\$^{b}$ Department of Mathematics
\\Eastern Michigan University
\\Ypsilanti, MI 48197, USA
\\
\\$^{c}$ Department of Mathematics
\\Northwestern Polytechnical University
\\Xi'an Shaanxi, 710072, PRC
\\
\\Email: jliu@emich.edu}
\date{ }
\maketitle

\dspace


\\

\begin{abstract}
In this paper, we provide a common generalization to the well-known Erd\H{o}s-Ko-Rado Theorem, Frankl-Wilson Theorem, Alon-Babai-Suzuki Theorem, and Snevily Theorem on set systems with $\mathcal{L}$-intersections. As a consequence, we derive a result which strengthens substantially the well-known theorem on set systems with $k$-wise $\mathcal{L}$-intersections by F$\ddot{u}$redi and Sudakov [J. Combin. Theory, Ser. A (2004) 105: 143-159]. We will also derive similar results on $\mathcal{L}$-intersecting families of subspaces of an $n$-dimensional vector space over a finite field $\mathbb{F}_{q}$, where $q$ is a prime power.
\end{abstract}

\section{Introduction}
A family $\mathcal{F}$ of subsets of $[n]=\{1, 2, \dots, n\}$ is called
$t$-$intersecting$ if for every pair of distinct subsets $E, F \in
\mathcal{F}$, $|E \cap F| \geq t$ (also called intersecting when $t = 1$). Let $\mathcal{L}=\{l_1,
l_2, \dots, l_s\}$ be a set of $s$ nonnegative integers.
A family $\mathcal{F}$ of subsets of $[n]$ is called $h$-$wise$
$\mathcal{L}$-$intersecting$ if $|F_{1} \cap F_{2} \cap \dots \cap F_{h}| \in \mathcal{L}$ for
every collection of $h$ distinct subsets in $\mathcal{F}$. When $h = 2$, such a family $\mathcal{F}$ is called $\mathcal{L}$-$intersecting$.
$\mathcal{F}$ is $k$-$uniform$ if it is a collection of $k$-subsets
of $[n]$. Thus, a $k$-uniform $t$-intersecting family is
$\mathcal{L}$-intersecting for $\mathcal{L}=\{t, t+1, \dots, k-1\}$.

In 1961, Erd\H{o}s, Ko, and Rado \cite{ekr} proved the following classical result.

\vspace{3mm}

\noindent {\bf Theorem 1.1} (Erd\H{o}s, Ko, and Rado, 1961 \cite{ekr}).  Let
$n \geq 2k$ and let $\mathcal{A}$ be a $k$-uniform intersecting
family of subsets of $[n]$. Then $|\mathcal{A}|\leq {{n-1} \choose
{k-1}}$ with equality only when $\mathcal{A}$ consists of all
$k$-subsets containing a common element.

\vspace{3mm}

To date, many intersection theorems have appeared in the literature, see \cite{ly} for a brief survey on theorems about $\mathcal{L}$-intersecting families.
The following $t$-intersecting version of Theorem 1.1 is due to Erd\H{o}s et al. \cite{ekr},  Frankl \cite{f}, and Wilson \cite{w}.

\vspace{3mm}

\noindent {\bf Theorem 1.2.}  Let
$n \geq (t+1)(k-t+1)$ and let $\mathcal{A}$ be a $k$-uniform $t$-intersecting
family of subsets of $[n]$. Then $|\mathcal{A}|\leq {{n-t} \choose
{k-t}}$ with equality only when $\mathcal{A}$ consists of all
$k$-subsets containing a common $t$-subset.

\vspace{3mm}

Here are well-known Frankl-Wilson theorem \cite{fw} and Alon-Babai-Suzuki theorem \cite{abs}.

\vspace{3mm}

\noindent {\bf Theorem 1.3} (Frankl and Wilson, 1981). Let
$\mathcal{L}=\{l_1, l_2, \dots, l_s\}$ be a set of $s$ nonnegative
integers. If $\mathcal{A}$ is an $\mathcal{L}$-intersecting family
of subsets of $[n]$, then
\[|\mathcal{A}|\leq {{n} \choose {s}}+ {{n}
\choose {s-1}}+ \cdots + {{n} \choose {0}}.\]

\vspace{3mm}

\noindent {\bf Theorem 1.4} (Alon, Babai, and Suzuki, 1991). Let
$\mathcal{L}=\{l_1, l_2, \dots, l_s\}$ be a set of $s$ nonnegative
integers and $K = \{k_1, k_2, \dots, k_r \}$ be a set of integers
satisfying $k_i > s-r$ for every $i$.
Suppose that $\mathcal{A} = \{A_1, A_2, \dots, A_{m}\}$ is a family of subsets of $[n]$ such that $|A_i| \in K$ for every $1 \leq i \leq m$
and $|A_{i} \cap A_{j}| \in \mathcal{L}$ for every pair $i \neq j$.
Then
\[m \leq {{n} \choose {s}}+ {{n} \choose {s-1}}+ \cdots + {{n} \choose {s-r+1}}.\]

\vspace{3mm}

Stronger bounds can be obtained if information about the specific set $\mathcal{L}$ is used. To that end,
Snevily \cite{s1} proved in 2003 the following theorem conjectured in 1994 by himself \cite{s}, which provides a common generalization of Frankl-F$\ddot{u}$redi theorem \cite{ff} (where $\mathcal{L}=\{1, 2, \dots, s\}$) and  Frankl-Wilson theorem (Theorem 1.3).

\vspace{3mm}

\noindent {\bf Theorem 1.5} (Snevily, 2003). Let
$\mathcal{L}=\{l_1, l_2, \dots, l_s\}$ be a set of $s$ positive
integers. If $\mathcal{A}$ is an $\mathcal{L}$-intersecting family
of subsets of $[n]$, then
\[|\mathcal{A}|\leq {{n-1} \choose {s}}+ {{n-1}
\choose {s-1}}+ \cdots + {{n-1} \choose {0}}.\]

\vspace{3mm}

The next conjecture proposed by J. Liu and X. Liu \cite{ll} provides a common generalization to all theorems above if it is true.

\vspace{3mm}

\noindent {\bf Conjecture 1.6.}
Let $\mathcal{L}=\{l_1, l_2, \dots, l_s\}$ be a set of $s$ nonnegative
integers with $l_{1} < l_{2} < \dots < l_{s}$ and $K = \{k_1, k_2, \dots, k_r \}$ be a set of positive integers satisfying $k_i > s-r$ for every $i$.
Suppose that $\mathcal{A} = \{A_1, A_2, \dots, A_{m}\}$ is a family of subsets of $[n]$ such that $|A_i| \in K$ for every $1 \leq i \leq m$
and $|A_{i} \cap A_{j}| \in \mathcal{L}$ for every pair $i \neq j$.
Then
\[m \leq {{n - l_{1}} \choose {s}}+ {{n - l_{1}} \choose {s-1}}+ \cdots + {{n - l_{1}} \choose {s-r+1}}.\]

\vspace{3mm}

The classical Erd\"{o}s-Ko-Rado theorem (Theorem 1.1) is the special case of Conjecture 1.6 with $l_{1} = 1$, $r = 1$, and $\mathcal{L}=\{1, 2, \dots, k-1\}$; Theorem 1.2 is the special case with $l_{1} = t$, $r = 1$, and $\mathcal{L}=\{t, t+1, \dots, k-1\}$; the Frankl-Wilson Theorem (Theorem 1.3) is the special case $l_{1} \geq 0$ and $r = n$; the well-known Alon-Babai-Suzuki theorem (Theorem 1.4) is the special case with $l_{1} \geq 0$; the Snevily's Theorem (Theorem 1.5) is the special case with $l_{1} \geq 1$ and $r = n$; and the well-known Ray-Chaudhuri-Wilson theorem \cite{rw} is the special case with $l_{1} \geq 0$ and $r = 1$. The bound in the conjecture is best possible as shown by the family of all subsets of $[n]$ with sizes at most $s+l_{1}$ and at least $s - r + 1 + l_{1}$ which contain all $1, 2, \dots, l_{1}$, where $\mathcal{L}=\{l_{1}, l_{1}+1, \dots, s + l_{1}-1\}$.

The main result by Heged\H{u}s \cite{h} recently shows that Conjecture 1.6 holds when $r = 1$ (i.e., for uniform families).  In this paper, we will prove the following theorem which shows that Conjecture 1.6 holds when $n \geq {{k^{2}} \choose {l_{1}+1}} + l_{1}$, where $k = \max\{k_{j}:1 \leq j \leq r\}$.

\vspace{3mm}

\noindent {\bf Theorem 1.7.}
Let $\mathcal{L}=\{l_1, l_2, \dots, l_s\}$ be a set of $s$ nonnegative
integers with $l_{1} < l_{2} < \dots < l_{s}$ and $K = \{k_1, k_2, \dots, k_r \}$ be a set of positive integers with $k_{i} > s - r + l_{1}$ for every $1 \leq i \leq r$.
Suppose that $\mathcal{A} = \{A_1, A_2, \dots, A_{m}\}$ is a family of subsets of $[n]$ such that $|A_i| \in K$ for every $1 \leq i \leq m$
and $|A_{i} \cap A_{j}| \in \mathcal{L}$ for every pair $i \neq j$. If $n \geq {{k^{2}} \choose {l_{1}+1}}s + l_{1}$,
where $k = \max\{|A_{j}|:1 \leq j \leq m\}$,
then
\[m \leq {{n - l_{1}} \choose {s}}+ {{n - l_{1}} \choose {s-1}}+ \cdots + {{n - l_{1}} \choose {s-r+1}}.\]
Moreover, if $n > {{k^{2}} \choose {l_{1}+1}}s + l_{1}$, the equality holds only if there exists an $l_{1}$-subset $T$ such that $T \subseteq A_{j}$ for every $j$.

\vspace{3mm}

As an immediate consequence of Theorem 1.7 by taking $r = n + 1$, we have the next corollary which improves Theorems 1.3 and 1.5 significantly when $n \geq {{k^{2}} \choose {l_{1}+1}}s + l_{1}$, where we have ${{n - l_{1}} \choose {i}} = 0$ if $i < 0$ and the obvious facts $s \leq n$ and $k_{i} \geq l_{1}$ for every $1 \leq i \leq r$.

\vspace{3mm}

\noindent {\bf Corollary 1.8.}
Let $\mathcal{L}=\{l_1, l_2, \dots, l_s\}$ be a set of $s$ nonnegative
integers with $l_{1} < l_{2} < \dots < l_{s}$.
Soppose that $\mathcal{A} = \{A_1, A_2, \dots, A_{m}\}$ is a family of subsets of $[n]$ such that
$|A_{i} \cap A_{j}| \in \mathcal{L}$ for every pair $i \neq j$. If $n \geq {{k^{2}} \choose {l_{1}+1}}s + l_{1}$,
where $k = \max\{|A_{j}|:1 \leq j \leq m\}$, then
\[m \leq {{n - l_{1}} \choose {s}}+ {{n - l_{1}} \choose {s-1}}+ \cdots + {{n - l_{1}} \choose {0}}.\]
Moreover, if $n > {{k^{2}} \choose {l_{1}+1}}s + l_{1}$, the equality holds only if there exists an $l_{1}$-subset $T$ such that $T \subseteq A_{j}$ for every $j$.

\vspace{3mm}

For set systems with $k$-wise $\mathcal{L}$-intersections,   F$\ddot{u}$redi and Sudakov \cite{fs} derived the following well-known theorem.

\vspace{3mm}

\noindent {\bf Theorem 1.9} (F$\ddot{u}$redi and Sudakov, 2004). Let $k \geq 3$ and let
$\mathcal{L}=\{l_1, l_2, \dots, l_s\}$ be a set of $s$ nonnegative
integers with $l_{1} < l_{2} < \dots < l_{s}$. If $\mathcal{A} = \{A_1, A_2, \dots, A_{m}\}$ is a family of subsets of $[n]$ such that
$|A_{i_1} \cap A_{i_2}\cap \dots \cap A_{i_k}| \in \mathcal{L}$ for every collection of $k$ distinct subsets in $\mathcal{A}$, then there exists $n_{0} = n_{0}(k, s)$ such that for all $n \geq n_{0}$
\[m \leq \frac{k + s - 1}{s + 1}{{n}\choose {s}}+ \sum_{ i \leq s-1}{{n} \choose {i}}.\]

\vspace{3mm}

As a consequence of Theorem 1.7 (Corollary 1.8), we will derive the next result which strengthens Theorem 1.9 substantially.

\vspace{3mm}

\noindent {\bf Theorem 1.10.} Let $k \geq 3$ and let
$\mathcal{L}=\{l_1, l_2, \dots, l_s\}$ be a set of $s$ nonnegative
integers with $l_{1} < l_{2} < \dots < l_{s}$. If $\mathcal{A} = \{A_1, A_2, \dots, A_{m}\}$ is a family of subsets of $[n]$ such that
$|A_{i_1} \cap A_{i_2}\cap \dots \cap A_{i_k}| \in \mathcal{L}$ for every collection of $k$ distinct subsets in $\mathcal{A}$, then there exists $n_{0} = n_{0}(k, s)$ such that for all $n \geq n_{0}$
\[m \leq \frac{k + s - 1}{s + 1}{{n-l_{1}}\choose {s}}+ \sum_{ i \leq s-1}{{n-l_{1}} \choose {i}}.\]

\vspace{3mm}

Next, we denote
\[
\qbinom{n}{k} := \prod_{0 \le i \leq k - 1} \frac{q^{n - i}-1}{q^{k - i}-1}.
\]
As for vector spaces over a finite field $\mathbb{F}_{q}$,
it is well-known that the number of all $k$-dimensional subspaces of an $n$-dimensional vector space over $\mathbb{F}_{q}$ is equal to $\qbinom{n}{k}$.

The next theorem is proved by Hsieh \cite{h1} for the case $n \ge 2k + 1$, $q \ge 3$ and $n \ge 2k + 2$, $q = 2$, and $\mathrm{dim}(V_{i} \cap V_{j}) > t - 1$; proved by Greene and Kleitman \cite{gk} for the case $n \ge 2k$, $k$ divides $n$; and proved completely by  Deza and Frankl \cite{df}.

\vspace{3mm}

\noindent {\bf Theorem 1.11} (Hsieh \cite{h1}, Greene and Kleitman \cite{gk}, Deza and Frankl \cite{df}).
Let $n \ge 2k$ and $\mathbb{F}_{q}$ be a finite field of order $q$. Suppose that $\mathcal{V}$ is a collection of $k$-subspaces of an $n$-dimensional vector space over $\mathbb{F}_{q}$ satisfying that $\mathrm{dim}(V_{i} \cap V_{j}) > 0$ for any distinct subspaces $V_{i}$ and $V_{j}$ in $\mathcal{V}$. Then
\[
|\mathcal{V}| \le \qbinom{n - 1}{k - 1}.
\]

\vspace{3mm}

In 1985, Frankl and Graham \cite{fg} proved the following theorem.

\vspace{3mm}

\noindent {\bf Theorem 1.12} (Frankl and Graham \cite{fg}). Let $\mathcal{L}$ be a set of $s$ nonnegative integers and $\mathbb{F}_{q}$ be a finite field of order $q$, $q$ is a prime power.
Suppose that $\mathcal{V}$ is a collection of $k$-dimensional subspaces of an $n$-dimensional vector space over $\mathbb{F}_{q}$ satisfying that $\mathrm{dim}(V_{i} \cap V_{j}) \in \mathcal{L}$ for any distinct subspaces $V_{i}$ and $V_{j}$ in $\mathcal{V}$. Then
\[
|\mathcal{V}| \le \qbinom{n}{s}.
\]

\vspace{3mm}

In 1990, Lefmann \cite{l} proved the following $\mathcal{L}$-intersecting theorem for ranked finite lattices. In 1993, Ray-Chaudhuri and Zhu \cite{rz} extended it to the polynomial semi-lattices. In 2001, Qian and Ray-Chaudhuri \cite{qr2} proved the following theorem for quasi-polynomial semi-lattices. The next result is a special case of their results.

\vspace{3mm}

\noindent {\bf Theorem 1.13} (Lefmann \cite{l}, Ray-Chaudhuri and Zhu \cite{rz}, Qian and Ray-Chaudhuri\cite{qr2}).
Let $\mathcal{L}$ be a set of $s$ nonnegative integers and $\mathbb{F}_{q}$ be a finite field of order $q$. Suppose that $\mathcal{V}$ is a collection of subspaces of an $n$-dimensional vector space over $\mathbb{F}_{q}$ satisfying that $\mathrm{dim}(V_{i} \cap V_{j}) \in \mathcal{L}$ for any distinct subspaces $V_{i}$ and $V_{j}$ in $\mathcal{V}$. Then
\[
|\mathcal{V}| \le \qbinom{n}{s} + \qbinom{n}{s - 1} + \cdots + \qbinom{n}{0}.
\]

\vspace{3mm}

Alon, Babai and Suzuki \cite{abs} derived the following stronger theorem in 1991.

\vspace{3mm}

\noindent {\bf Theorem 1.14} (Alon, Babai and Suzuki \cite{abs}).
Let $\mathcal{L}$ be a set of $s$ nonnegative integers and $\mathbb{F}_{q}$ be a finite field of order $q$. Suppose that $\mathcal{V}$ is a collection of subspaces of an $n$-dimensional vector space over $\mathbb{F}_{q}$ satisfying that $\mathrm{dim}(V_{i} \cap V_{j}) \in \mathcal{L}$ for any distinct subspaces $V_{i}$ and $V_{j}$ in $\mathcal{V}$ and the dimension of every member of $\mathcal{V}$ belongs to the set $\{k_{1}, k_{2}, \ldots, k_{t}\}$ with $k_{i} > s - t$ for every $i$. Then
\[
|\mathcal{V}| \le \qbinom{n}{s} + \qbinom{n}{s - 1} + \cdots + \qbinom{n}{s - t + 1}.
\]

\vspace{3mm}

Here we will prove the next two results which strengthen Theorems 1.13 and 1.14.

\vspace{3mm}

\noindent {\bf Theorem 1.15.}
Let $\mathcal{L}$ be a set of $s$ nonnegative integers with $l_{1} < l_{2} < \cdots <l_{s}$ and $\mathbb{F}_{q}$ be a finite field of order $q$. Suppose that $\mathcal{V} = \{V_{1}, V_{2}, \ldots, V_{m}\}$ is a collection of subspaces of an $n$-dimensional vector space over $\mathbb{F}_{q}$ satisfying that $\mathrm{dim}(V_{i} \cap V_{j}) \in \mathcal{L}$ for any distinct subspaces $V_{i}$ and $V_{j}$ in $\mathcal{V}$.  If $n \ge \log_{q}((q^{s} - 1)\qbinom{k^{2}}{l_{1} + 1} + 1) + l_{1}$, where $k = \max\{\dim(V_{j}) : 1 \le j \le m\}$. Then
\[
|\mathcal{V}| \le \qbinom{n - l_{1}}{s} + \qbinom{n - l_{1}}{s - 1} + \cdots + \qbinom{n - l_{1}}{0}.
\]

\vspace{3mm}

\noindent {\bf Theorem 1.16.}
Let $\mathcal{L}$ be a set of $s$ nonnegative integers with $l_{1} < l_{2} < \cdots <l_{s}$ and $\mathbb{F}_{q}$ be a finite field of order $q$. Suppose that $\mathcal{V} = \{V_{1}, V_{2}, \ldots, V_{m}\}$ is a collection of subspaces of an $n$-dimensional vector space over $\mathbb{F}_{q}$ satisfying that $\mathrm{dim}(V_{i} \cap V_{j}) \in \mathcal{L}$ for any distinct subspaces $V_{i}$ and $V_{j}$ in $\mathcal{V}$ and the dimension of every member of $\mathcal{V}$ belongs to the set $\{k_{1}, k_{2}, \ldots, k_{r}\}$ and $k_{i} > s - r + l_{1}$ for every $i$. If $n \ge \log_{q}((q^{s} - 1)\qbinom{k^{2}}{l_{1} + 1} + 1) + l_{1}$, , where $k = \max\{\dim(V_{j}) : 1 \le j \le m\}$. Then
\[
|\mathcal{V}| \le \qbinom{n - l_{1}}{s} + \qbinom{n - l_{1}}{s - 1} + \cdots + \qbinom{n - l_{1}}{s - r + 1}.
\]

\vspace{3mm}

This paper is motivated by ideas  from \cite{fs} and \cite{h}.

\section{Proof of Theorem 1.7}
Throughout this paper, we use ${{[n]}\choose {k}}$ to denote the set of all $k$-subsets of $[n] = \{1, 2, \dots, n\}$.
A family $\mathcal{F}$ od sets is said to be an $H_{k}$-family ($k \geq 1$) if for every $\mathcal{G} \subseteq \mathcal{F}$, $\cap_{G \in \mathcal{G}}G = \emptyset$ implies that $\cap_{G \in \mathcal{G}'}G = \emptyset$ for some $\mathcal{G}' \subseteq \mathcal{G}$ with $|\mathcal{G}'| \leq k$. The following lemma is Theorem 1(i) from \cite{bd}.

\vspace{3mm}

\noindent {\bf Lemma 2.1.} If $\mathcal{F} \subseteq \cup_{i = 0}^{k}{{[n]}\choose {i}}$, then $\mathcal{F}$ is an $H_{d}$-family  for every $d \geq k + 1$, that is, if $\cap_{F \in \mathcal{F}}F = \emptyset$, then $\cap_{F \in \mathcal{F}'}F = \emptyset$ for some $\mathcal{F}' \subseteq \mathcal{F}$ with $|\mathcal{F}'| \leq k+1$.

\vspace{3mm}

The next two Lemmas are Lemmas 2.2 and 2.3 from \cite{h}.

\vspace{3mm}

\noindent {\bf Lemma 2.2.} Let $l_{1}$ be a positive integer. Let $\mathcal{H}$ be a family of subsets of $[n]$. Suppose that $\cap_{H \in \mathcal{H}}H = \emptyset$. Let $F \subseteq [n]$, $F \not \in \mathcal{H}$ be a subset such that $|F \cap H| \geq l_{1}$ for each $H \in \mathcal{H}$. Set $Q = \cup_{H \in \mathcal{H}}H$. Then
\[|Q \cap F| \geq l_{1} + 1.\]

\vspace{3mm}

\noindent {\bf Lemma 2.3.} Let $\mathcal{H}$ be a family of subsets of $[n]$. Suppose that $t = |\mathcal{H}| \geq 2$ and $\mathcal{H}$ is a $k$-uniform intersecting family. Then
\[|\cup_{h \in \mathcal{H}}H| \leq k + (t - 1)(k - 1).\]

\vspace{3mm}

When $\mathcal{H}$ is a non-uniform family of subsets of $[n]$ with maximum sbuset size $k$, we can obtain the next corollary from Lemma 2.3 by  extending every subset to a $k$-subset in an arbitrary way.

\vspace{3mm}

\noindent {\bf Corollary 2.4.} Let $\mathcal{H}$ be a family of subsets of $[n]$ with maximum subset size $k$. Suppose that $t = |\mathcal{H}| \geq 2$ and $\mathcal{H}$ is an intersecting family. Then
\[|\cup_{h \in \mathcal{H}}H| \leq k + (t - 1)(k - 1).\]

\vspace{3mm}

\noindent {\bf Proof of Theorem 1.7.} Denote $k = \max\{|A_{j}|:1 \leq j \leq m\}$. By Theorem 1.4, we may assume that $l_{1} > 0$. Also, since $\mathcal{A}$ is $\mathcal{L}$-intersecting, if there exists $A_{i} \in \mathcal{A}$ such that $|A_{i}| = l_{1}$, then $A_{i} \subseteq \cap_{A_{j} \in \mathcal{A}}A_{j}$ and the result follows from Theorem 1.4 easily by considering the family $\{A_{j} \setminus A_{i}: 1 \leq j \leq m\}$ on the set $[n] \setminus A_{i}$, where $k_{j} -l_{1} > s - r$ for every $1 \leq j \leq r$. Thus, we assume that $|A_{j}| \geq l_{1}+1$ for all $A_{j} \in \mathcal{A}$.

We consider the following cases:

\noindent {\bf Case 1.} $\cap_{A_{i} \in \mathcal{A}}A_{i} = \emptyset$. By Lemma 2.1, there exists a subfamily $\mathcal{F} \subseteq \mathcal{A}$ with $|\mathcal{F}| = k+1$ such that $\cap_{A_{j} \in \mathcal{F}}A_{j} = \emptyset$. Let
\[M = \cup_{A_{j} \in \mathcal{F}}A_{j}.\]
Then $|M| \leq k + k(k-1) = k^{2}$ by Corollary 2.4. On the other hand, since $|A_{i}| \geq l_{1}+1$ for all $A_{i} \in \mathcal{A}$, it follows from Lemma 2.2 that
\[\hspace{39mm} |M \cap A_{i}| \geq l_{1} + 1 \mbox{ for each } A_{i} \in \mathcal{A}. \hspace{39mm} (2.1)\]
Let $T$ be a given subset of $M$ such that $|T| = l_{1}+1$. Define
\[\mathcal{A}(T) = \{A_{i} \in \mathcal{A}: T \subseteq M \cap A_{i}\}.\]
Set $\mathcal{L}'=\{l_2, l_3, \dots, l_s\}$. Then $|\mathcal{L}'| = s - 1$. Since  $\mathcal{A}$ is $\mathcal{L}$-intersecting family and $|E \cap F| \geq |T| \geq l_{1}+1$ for any $E, F \in \mathcal{A}(T)$, $\mathcal{A}(T)$ is  $\mathcal{L}'$-intersecting.
By (2.1), it is easy to check that
\[\hspace{48mm} \mathcal{A} = \cup_{T \subseteq M, |T| = l_{1}+1}\mathcal{A}(T). \hspace{48mm} (2.2)\]
Note that for each $T \subseteq M$ with $|T| = l_{1}+1$, the system
\[\mathcal{G}(T) = \{A_{i} \setminus T: A_{i} \in \mathcal{A}(T) \}\]
is an $\mathcal{L}^{*}$-intersecting family on the set $[n] \setminus T$, where $\mathcal{L}^{*} = \{l_2 - l_{1} - 1, l_3  - l_{1} - 1, \dots, l_s  - l_{1} - 1\}$, and $|\mathcal{G}(T)| = |\mathcal{A}(T)|$. Since $k_{j} - (l_{1} + 1) > s - 1 - r$ for each $1 \leq j \leq r$, it follows from Theorem 1.4 that for each $T \subseteq M$ with $|T| = l_{1}+1$,
\[|\mathcal{A}(T)| \leq \sum_{j = s - 1 - r + 1}^{s - 1}{{n - l_{1} - 1}\choose {j}}.\]
Since $n \geq {{k^{2}} \choose {l_{1}+1}}s + l_{1}$, it follows from (2.2) that
\[|\mathcal{A}| \leq \sum_{T \subseteq M, |T| = l_{1}+1}|\mathcal{A}(T)| \leq {{k^{2}}\choose {l_{1}+1}}\sum_{j = s - 1 - r + 1}^{s - 1}{{n - l_{1} - 1}\choose {j}} \]
\[\leq {{k^{2}}\choose {l_{1}+1}}\frac{s}{n - l_{1}}\sum_{j = s - 1 - r + 1}^{s - 1}{{n - l_{1}}\choose {j+1}} \leq \sum_{j = s - r + 1}^{s}{{n - l_{1}}\choose {j}}.\]

\noindent {\bf Case 2.} $\cap_{A_{i} \in \mathcal{A}}A_{i} \neq \emptyset$. For this case, if $|\cap_{A_{i} \in \mathcal{A}}A_{i}| \geq l_{1}$, then the result follows easily from Theorem 1.4. Assume that $0 < |\cap_{A_{i} \in \mathcal{A}}A_{i}| = t < l_{1}$ and let $T = \cap_{A_{i} \in \mathcal{A}}A_{i}$. Then
\[\mathcal{G} = \{A_{i} \setminus T: A_{i} \in \mathcal{A} \}\]
is an $\mathcal{L}'$-intersecting family on the set $[n] \setminus T$, where $\mathcal{L}'=\{l_1 - t, l_2 - t, \dots, l_s - t\}$ with $l_{1} - t > 0$. By Case 1, we obtain
\[|\mathcal{A}| = |\mathcal{G}| \leq \sum_{j = s - r + 1}^{s}{{(n - t) - (l_{1} - t)}\choose {j}} = \sum_{j = s - r + 1}^{s}{{n - l_{1}}\choose {j}}.\]
From the arguments above, it is clear that if $n > {{k^{2}} \choose {l_{1}+1}}s + l_{1}$, the equality holds only if there exists an $l_{1}$-subset $T$ such that $T \subseteq A_{j}$ for every $j$.
\hfill$\Box$

\section{Proof of Theorem 1.10}
We begin with the following lemma which follows easily from Lemma 2.1.

\vspace{3mm}

\noindent {\bf Lemma 3.1.} Let $\mathcal{F} \subseteq \cup_{i = 0}^{k}{{[n]}\choose {i}}$.  If $\cap_{F \in \mathcal{F}}F = T$, then $\cap_{F \in \mathcal{F}'}F = T$ for some $\mathcal{F}' \subseteq \mathcal{F}$ with $|\mathcal{F}'| \leq k+1$.

\vspace{3mm}

The next lemma is Lemma 3.1 in \cite{fs}.

\vspace{3mm}

\noindent {\bf Lemma 3.2.}
Let $\mathcal{L}=\{l_1, l_2, \dots, l_s\}$ be a set of $s$ nonnegative integers.
Suppose that $\mathcal{A} = \{A_1, A_2, \dots, A_{m}\}$ and
$\mathcal{B} = \{B_{1}, B_{2}, \dots, B_{m}\}$ are two families of subsets of $[n]$ such that
(i) $|A_{i} \cap B_{j}|  \in \mathcal{L}$ for every pair $i < j$,
(ii)  $|A_{i} \cap B_{i}| \not \in \mathcal{L}$ for any $1 \leq i \leq m$.
Then
\[m \leq {{n}\choose {s}}+ {{n} \choose {s-1}}+ \cdots + {{n} \choose {0}}.\]

\vspace{3mm}

\vspace{3mm}

To prove Theorem 1.10, we need the following result.

\vspace{3mm}

\noindent {\bf Proposition 3.3.}
Let $\mathcal{L}=\{l_1, l_2, \dots, l_s\}$ be a set of $s$
integers with $0 < l_{1} < l_{2} < \dots < l_{s}$.
Suppose that $\mathcal{A} = \{A_1, A_2, \dots, A_{m}\}$ and
$\mathcal{B} = \{B_{1}, B_{2}, \dots, B_{m}\}$ are two families of subsets of $[n]$ such that
(i) $|A_{i} \cap B_{j}|  \in \mathcal{L}$ for every pair $i < j$,
(ii)  $B_{i} \subseteq A_{i}$ for all $i \leq m$ and $|A_{i} \cap B_{i}| \not \in \mathcal{L}$ for any $k + 2 \leq i \leq m$,
where $k$ is the maximum subset size in $\mathcal{B}$,
(iii) $|\cap_{1 \leq j \leq k+1} B_{j}| = |\cap_{B_j \in \mathcal{B}} B_{j}| < l_{1}$  and  $A_{i} = B_{i}$ for every $i \leq k+1$.
If $n \geq [{{k^{2} + k} \choose {l_{1}+1}}+1]s + l_{1}$,  then
\[m \leq {{n-l_{1}}\choose {s}}+ {{n-l_{1}} \choose {s-1}}+ \cdots + {{n-l_{1}} \choose {0}}.\]

\vspace{3mm}

\noindent {\bf Proof.}
We consider the following cases:

\noindent {\bf Case 1.} $\cap_{B_j \in \mathcal{B}}B_{j} = \emptyset$. By assumption (iii), we have  $\cap_{j \leq k+1}B_{j} = \emptyset$, where  $k$ is the maximum subset size in $\mathcal{B}$. Let
\[M = \cup_{j \leq k+1}B_{j} = \cup_{j \leq k+1}A_{j}.\]
Then $|M| \leq k(k+1) = k^{2} + k$. On the other hand, it follows from Lemma 2.2 and assumptions (i) and (iii) that
\[\hspace{26mm} |M \cap B_{j}| \geq l_{1} + 1 \mbox{ for each } B_{j} \in \mathcal{B} \mbox{ with } j \geq k + 2. \hspace{26mm} (3.1)\]
Let $T$ be a given subset of $M$ such that $|T| = l_{1}+1$. Define
\[\mathcal{B}(T) = \{B_{j} \in \mathcal{B}: j \geq k+2 \mbox{ and } T \subseteq M \cap B_{j}\} \]
and set
\[\mathcal{A}(T) = \{A_{j} \in \mathcal{A}: j \geq k+2 \mbox{ and } B_{j} \in \mathcal{B}(T)\}. \]
Set $\mathcal{L}'=\{l_2, l_3, \dots, l_s\}$. Then $|\mathcal{L}'| = s - 1$ and $\mathcal{A}(T)$ and $\mathcal{B}(T)$ satisfy (i) $|A_{i} \cap B_{j}|  \in \mathcal{L}'$ for every pair $i < j$, (ii) $T \subseteq B_{j} \subseteq A_{j}$ and $|A_{j} \cap B_{j}| \not \in \mathcal{L}$ for each $k+2 \leq j \leq m$ such that $B_{j} \in \mathcal{B}(T)$. Define
\[\mathcal{G}(T) = \{A_{j} \setminus T: A_{j} \in \mathcal{A}(T) \},\]
\[\mathcal{Q}(T) = \{B_{j} \setminus T: B_{j} \in \mathcal{B}(T) \}.\]
Then $\mathcal{G}(T)$ and $\mathcal{Q}(T)$ are two families on $[n] \setminus T$ satisfying
that $|(A_{i}\setminus T) \cap (B_{j} \setminus T)| \in \mathcal{L}^{*} = \{l_2 - l_{1} - 1, l_3  - l_{1} - 1, \dots, l_s  - l_{1} - 1\}$ for $A_{i} \in \mathcal{A}(T)$ and $B_{j} \in \mathcal{B}(T)$ with $i < j$.
It follows from Lemma 3.2 that
\[|\mathcal{B}(T)| = |\mathcal{Q}(T)| \leq \sum_{j = 0}^{s - 1}{{n - l_{1} - 1}\choose {j}}.\]
By (3.1), it is easy to check that
\[\hspace{33mm} \mathcal{B}\setminus \{B_{i}: i \leq k+1\} =  \cup_{T \subseteq M, |T| = l_{1}+1}\mathcal{B}(T). \hspace{33mm} (3.2)\]
Since $n \geq [{{k^{2}} \choose {l_{1}+1}}+1]s + l_{1}$, it follows from (3.2) that
\[m = |\mathcal{B}| \leq \sum_{T \subseteq M, |T| = l_{1}+1}|\mathcal{B}(T)| + k + 1
\leq {{k^{2}}\choose {l_{1}+1}}\sum_{j = 0}^{s - 1}{{n - l_{1} - 1}\choose {j}} + k + 1 \]
\[\leq \bigg[{{k^{2}}\choose {l_{1}+1}}+1\bigg]\sum_{j = 0}^{s - 1}{{n - l_{1} - 1}\choose {j}}\]
\[\leq \bigg[{{k^{2}}\choose {l_{1}+1}}+1\bigg]\frac{s}{n - l_{1}}\sum_{j = 0}^{s - 1}{{n - l_{1}}\choose {j+1}}  \]
\[\leq \sum_{j = 1}^{s}{{n - l_{1}}\choose {j}}  \leq \sum_{j = 0}^{s}{{n - l_{1}}\choose {j}}.\]

\noindent {\bf Case 2.} $\cap_{B_{j} \in \mathcal{B}}B_{j} \neq \emptyset$. Let $T = \cap_{B_{j} \in \mathcal{B}}B_{j}$ and $t = |T|$.
By assumption (iii),  $0 < t < l_{1}$ and  $\cap_{j \leq k+1}B_{j} = T$.
Set
\[\mathcal{G} = \{A_{i} \setminus T: A_{i} \in \mathcal{A} \}\]
\[\mathcal{H} = \{B_{i} \setminus T: B_{i} \in \mathcal{B} \}.\]
Then $\mathcal{G}$ and $\mathcal{H}$ are two families on $[n] \setminus T$ satisfying that
for any $G_{i} = A_{i} \setminus T \in \mathcal{G}$ and $H_{j} = B_{j} \setminus T \in \mathcal{H}$ with $i < j$, $|G_{i} \cap H_{j}| \in \mathcal{L}'$, where $\mathcal{L}'=\{l_1 - t, l_2 - t, \dots, l_s - t\}$ with $l_{1} - t > 0$.
Moreover, $\mathcal{G}$ and $\mathcal{H}$ satisfy assumptions (ii) and (iii).
By Case 1, we obtain
\[|\mathcal{A}| = |\mathcal{G}| \leq \sum_{j = 0}^{s}{{(n - t) - (l_{1} - t)}\choose {j}} = \sum_{j = 0}^{s}{{n - l_{1}}\choose {j}}.\]
\hfill$\Box$

\vspace{3mm}

We also need the next result by Grolmusz and Sudakov \cite{gs} which extends Theorem 1.3
to $k$-wise $\mathcal{L}$-intersecting families.

\vspace{3mm}

\noindent {\bf Theorem 3.4} (Grolmusz and Sudakov, 2002 \cite{gs}). Let $k \geq 2$ and let
$\mathcal{L}=\{l_1, l_2, \dots, l_s\}$ be a set of $s$ nonnegative
integers. If $\mathcal{A} = \{A_1, A_2, \dots, A_{m}\}$ is a family
of subsets of $[n]$ such that $|A_{i_1} \cap A_{i_2}\cap \dots \cap A_{i_k}| \in \mathcal{L}$ for every collection of $k$ distinct subsets in $\mathcal{A}$, then
\[m \leq (k - 1)\Bigg[{{n}\choose {s}}+ {{n} \choose {s-1}} + \cdots + {{n} \choose {0}}\Bigg].\]

\vspace{3mm}

The following theorem strengthens Theorem 3.4 considerably.

\vspace{3mm}

\noindent {\bf Theorem 3.5.} Let $k \geq 2$,
$\mathcal{L}=\{l_1, l_2, \dots, l_s\}$ be a set of $s$ nonnegative integers with $l_{1} < l_{2} < \dots < l_{s}$.
Suppose that $\mathcal{A} = \{A_1, A_2, \dots, A_{m}\}$ is a family of subsets of $[n]$ such that  $|A_{i_1} \cap A_{i_2}\cap \dots \cap A_{i_k}| \in \mathcal{L}$ for every collection of $k$ distinct subsets in $\mathcal{A}$.
If $n \geq [{{h^{2} + h} \choose {l_{1}+1}}+1]s + l_{1}$, where $h$ is the maximum subset size in $\mathcal{A}$, then
\[m \leq (k - 1)\Bigg[{{n-l_{1}}\choose {s}}+ {{n-l_{1}} \choose {s-1}}+
\cdots + {{n-l_{1}} \choose {0}}\Bigg].\]

\vspace{3mm}

\noindent {\bf Proof.} The case $l_{1} = 0$ follows from Theorem 3.4. So we assume that $l_{1} > 0$. If $|\cap_{A_{j} \in \mathcal{A}} A_{j}| \geq l_{1}$, then the result also follows easily from Theorem 3.4 by considering the family $\mathcal{A}' = \{ A_{j} \setminus T : A_{j} \in \mathcal{A}\}$, where $T \subseteq \cap_{A_{j} \in \mathcal{A}} A_{j}$ with $|T| = l_{1}$. Thus, we assume that $|\cap_{A_{j} \in \mathcal{A}} A_{j}| < l_{1}$.

We now proceed by induction on $k \geq 2$. The case $k = 2$ is just Corollary 1.8. Assume that the theorem holds for $(k-1)$-wise $\mathcal{L}$-intersecting families. Next, we show that the theorem holds for $k$-wise $\mathcal{L}$-intersecting families.
We need to partition $\mathcal{A}$ into two families $\mathcal{B}$ and $\mathcal{F}$ with the following properties: $\mathcal{F}$ is $(k-1)$-wise $\mathcal{L}$-intersecting and there exists a family $\mathcal{C}$ such that the pair $(\mathcal{B}, \mathcal{C})$ satisfies the assumptions in Proposition 3.3.

To obtain the desired partition, we first construct families $\mathcal{F}$, $\mathcal{B}$, and $\mathcal{C}$ by repeating the following procedure:
By Lemma 3.1, without loss of generality,  we have $\cap_{j \leq k + 1}A_{j} = \cap_{A_{j} \in \mathcal{A}} A_{j}$(recall that $|\cap_{A_{j} \in \mathcal{A}} A_{j}| < l_{1}$). Set $B_{i} = C_{i} = A_{i}$ for $i \leq k + 1$.
For every $k+1 \leq d \leq m-1$, suppose that after step $d$ we have constructed families $\mathcal{B} = \{B_{1}, B_{2}, \dots, B_{d}\}$, $\mathcal{C} = \{C_{1}, C_{2}, \dots, C_{d}\}$. At step $d+1$:
suppose there are indices $i_{1} < i_{2} < \cdots < i_{k-1}$ such that $|A_{i_1} \cap A_{i_2}\cap \dots \cap A_{i_{k-1}}| \not \in \mathcal{L}$. With relabeling if necessary, we assume $i_{1} = d+1$
Then set $B_{d+1} = A_{i_{1}} = A_{d+1}$, $C_{d+1} = A_{i_1} \cap A_{i_2}\cap \dots \cap A_{i_{k-1}}$.  Note that $C_{d+1} \subseteq B_{d+1}$ and  $|B_{d+1} \cap C_{d+1}| = |C_{d+1}| \not \in \mathcal{L}$ and $|B_{j} \cap C_{d+1}| \in  \mathcal{L}$ for all $j < d+1$. Update $d$ by $d+1$ and proceed to the next step. Continue this process until we can not proceed further. Then set $\mathcal{F} = \mathcal{A} \setminus \mathcal{B}$.

Then, by the construction, $m \leq |\mathcal{B}|+ |\mathcal{F}| = m' + |\mathcal{F}|$, the families $\mathcal{B}$ and $\mathcal{C}$ satisfy the assumptions (i) - (iii) in Proposition 3.3 (with $\mathcal{B}$ replacing $\mathcal{A}$ and $\mathcal{C}$ replacing $\mathcal{B}$), and $\mathcal{F}$ is $(k-1)$-wise $\mathcal{L}$-intersecting. It follows from Proposition 3.3 that
\[m' \leq {{n-l_{1}}\choose {s}}+ {{n-l_{1}} \choose {s-1}} + \cdots + {{n-l_{1}} \choose {0}}.\]
By the induction hypothesis,
\[|\mathcal{F}| \leq (k-2)\Bigg[{{n-l_{1}}\choose {s}}+ {{n-l_{1}} \choose {s-1}} + \cdots + {{n-l_{1}} \choose {0}}\Bigg].\]
It follows that
\[m = m' + |\mathcal{F}| \leq (k-1)\Bigg[{{n-l_{1}}\choose {s}}+ {{n-l_{1}} \choose {s-1}} + \cdots + {{n-l_{1}} \choose {0}}\Bigg]. \]
\hfill$\Box$

\vspace{3mm}

Similar to the proof for Proposition 4.1 in F$\ddot{u}$redi and  Sudakov \cite{fs},  we can derive the next lemma by using Theorem 3.5 and Proposition 3.3.

\vspace{3mm}

\noindent {\bf Lemma 3.6.} Let $k \geq 3$ and $\mathcal{L}=\{l_1, l_2, \dots, l_s\}$ be a set of $s$ integers satisfying $0 < l_{1} < l_{2} < \dots < l_{s}$.
Suppose that $\mathcal{A} = \{A_1, A_2, \dots, A_{m}\}$ is a $k$-wise $\mathcal{L}$-intersecting family of subsets of $[n]$ and $n \geq {{h^{2} + h} \choose {l_{1}+1}}s + l_{1}$, where $h$ is the maximum subset size in $\mathcal{A}$. If there exists an index $r$, $1 \leq r \leq s$, such that no intersection of $k-1$ distinct members of $\mathcal{A}$ has size $l_{r}$, then
\[m \leq {{n-l_{1}}\choose {s}} + (k - 1)\sum_{i \leq s-1}{{n-l_{1}}\choose {i}}.\]

\noindent {\bf Proof.} Let $k \geq 3$ and $\mathcal{L}=\{l_1, l_2, \dots, l_s\}$ be a set of $s$ positive integers $l_{1} < l_{2} < \dots < l_{s}$
and let $\mathcal{A} = \{A_1, A_2, \dots, A_{m}\}$ be a $k$-wise $\mathcal{L}$-intersecting family of subsets of $[n]$.
Similar to the argument in the proof of Theorem 3.5, we partition $\mathcal{A}$ into two families $\mathcal{A}'$ and $\mathcal{F}$ with the following properties: $\mathcal{F}$ is $(k-1)$-wise $\mathcal{L}'$-intersecting with $\mathcal{L}' = \mathcal{L} \setminus \{l_{r}\}$ and there exists a family $\mathcal{B}'$ such that the pair $(\mathcal{A}', \mathcal{B}')$ satisfies the conditions of Proposition 3.3.
Applying Theorem 3.5 to the family $\mathcal{F}$, we obtain
\[|\mathcal{F}| \leq (k - 2)\Bigg[{{n-l_{1}}\choose {s-1}}+ {{n-l_{1}} \choose {s-2}} + \cdots + {{n-l_{1}} \choose {0}}\Bigg].\]
By Proposition 3.3, we have
\[m' = |\mathcal{A}'| \leq {{n-l_{1}}\choose {s}}+ {{n-l_{1}} \choose {s-1}}+ \cdots + {{n-l_{1}} \choose {0}}.\]
It follows that
\[|\mathcal{A}| = |\mathcal{F}| + |\mathcal{A}'| \leq {{n-l_{1}}\choose {s}} + (k - 1)\sum_{i \leq s-1}{{n-l_{1}}\choose {i}}.\]
\hfill$\Box$

\vspace{3mm}

As a consequence to Lemma 3.6, we have the following corollary.

\vspace{3mm}

\noindent {\bf Corollary 3.7.} Let $k \geq 3$ and $\mathcal{L}=\{l_1 < l_2 < \cdots < l_s\}$ be a set of $s$ positive integers
and let $\mathcal{A} = \{A_1, A_2, \dots, A_{m}\}$ be a $k$-wise $\mathcal{L}$-intersecting family of subsets of $[n]$. If
\[|\mathcal{A}| > {{n-l_{1}}\choose {s}} + (k - 1)\sum_{i \leq s-1}{{n-l_{1}}\choose {i}},\]
then there exists an $l_{1}$-set $X$ of $[n]$ such that $X \subseteq A_{j}$ for all $A_{j} \in \mathcal{A}$.

\vspace{3mm}

\noindent {\bf Proof.} Let $\mathcal{L}=\{l_1 < l_2 < \cdots < l_s\}$ be such that $ l_1 > 0$. If no intersection of $k-1$ distinct members of $\mathcal{A}$ has size $l_{1}$, then, by Lemma 3.6, we have $|\mathcal{A}| \leq {{n-l_{1}}\choose {s}} + (k - 1)\sum_{i \leq s-1}{{n-l_{1}}\choose {i}}$,
a contradiction. Thus, there are $k-1$ distinct members in $\mathcal{A}$, say $A_{1}$, $A_{2}$, \dots, $A_{k-1}$, such that $|A_{1} \cap A_{2} \cap \cdots \cap A_{k-1}| = l_{1}$. Since $\mathcal{A}$ is $k$-wise $\mathcal{L}$-intersecting and $l_1 < l_2 < \cdots < l_s$, it follows that
\[|A \cap (A_{1} \cap A_{2} \cap \cdots \cap A_{k-1})| = l_{1} \mbox{ for any other } A \in \mathcal{A}\]
which implies that $X = A_{1} \cap A_{2} \cap \cdots \cap A_{k-1} \subseteq A_{j}$ for all members $A_{j} \in \mathcal{A}$.
\hfill$\Box$

\vspace{3mm}

We are now ready to prove Theorem 1.10.

\vspace{3mm}

\noindent {\bf Proof of Theorem 1.10.}
Let $k \geq 3$ and $\mathcal{L}=\{l_1, l_2, \dots, l_s\}$ be a set of $s$ nonnegative integers with $l_1 < l_2 < \cdots < l_s$
and let $\mathcal{A} = \{A_1, A_2, \dots, A_{m}\}$ be a $k$-wise $\mathcal{L}$-intersecting family of subsets of $[n]$.
By Theorem 1.9, we may assume that $l_{1} > 0$.
Note that for $n$ sufficiently large,
\[{{n-l_{1}}\choose {s}} + (k - 1)\sum_{i \leq s-1}{{n-l_{1}}\choose {i}} \leq \frac{k + s - 1}{s + 1}{{n-l_{1}}\choose {s}}+ \sum_{i \leq s-1}{{n-l_{1}} \choose {i}}.\]
If $|\mathcal{A}| \leq  {{n-l_{1}}\choose {s}} + (k - 1)\sum_{i \leq s-1}{{n-l_{1}}\choose {i}}$, then the theorem follows.
Suppose that
\[|\mathcal{A}| > {{n-l_{1}}\choose {s}} + (k - 1)\sum_{i \leq s-1}{{n-l_{1}}\choose {i}}.\]
By Corollary 3.7, there exists an $l_{1}$-set $X$ of $[n]$ such that $X \subseteq A_{j}$ for all $A_{j} \in \mathcal{A}$. Then the family
$\mathcal{A}' = \{A_{j}\setminus X : A_{j} \in \mathcal{A}\}$ is a $k$-wise $\mathcal{L}'$-intersecting family of subsets of an $(n-l_{1})$-element set with $|\mathcal{A}| = |\mathcal{A}'|$ and $\mathcal{L}' = \{0, l_2 - l_1, \dots, l_{s} - l_{1}\}$. It follows from Theorem 1.9 that for $n \geq n_{0}$,
\[|\mathcal{A}| = |\mathcal{A}'| \leq \frac{k + s - 1}{s + 1}{{n-l_{1}}\choose {s}}+ \sum_{i \leq s-1}{{n-l_{1}} \choose {i}}.\]
\hfill$\Box$

\vspace{3mm}

\section{Proof of Theorems 1.15 and 1.16}
We write $U \subseteq V$ if $U$ is a subspace of $V$ and denote the $n$-dimensional vector space by $W$.

\vspace{3mm}

\noindent {\bf Lemma 4.1.}
Let $\mathbb{F}_{q}$ be a finite field of order $q$ and $\mathcal{V}$ be a collection of subspaces of an $n$-dimensional vector space over $\mathbb{F}_{q}$ such that the dimension of every member of $\mathcal{V}$ is at most $k$. If $\dim(\cap_{V\in\mathcal{V}} V) = 0$, then $\dim(\cap_{V\in\mathcal{V^{\prime}}} V) = 0$ for some $\mathcal{V}^{\prime} \subseteq \mathcal{V}$ with $|\mathcal{V}^{\prime}| \le k + 1.$

\vspace{3mm}

\noindent {\bf Proof.} Suppose that $\mathcal{V}$ is a collection of subspaces of an $n$-dimensional vector space $W$ over $\mathbb{F}_{q}$ such that the dimension of every member of $\mathcal{V}$ is at most $k$ and $\dim(\cap_{V\in\mathcal{V}} V) = 0$. Let $\mathcal{V}^{\prime}$ be a subfamily of $\mathcal{V}$ of minimum size which satisfies $\dim(\cap_{V\in\mathcal{V^{\prime}}} V) = 0$. Assume that  $|\mathcal{V}^{\prime}| = p + 1$ with $p \geq 0$ and $\mathcal{V}^{\prime} = \{V_{1}, V_{2}, \dots, V_{p+1}\}$. We now show that $p \leq k$. Clearly, if $p = 0$, then $p \leq k$. So, we assume that $p \geq 1$. Denote  $\mathcal{V}_{i} = \mathcal{V}^{\prime} \setminus \{V_{i}\}$ for $1 \leq i \leq p + 1$. Then $\dim(\cap_{V\in\mathcal{V}_{i}} V) \not= 0$ for $i = 1, 2, \ldots, p+1$ by the choice of $\mathcal{V}^{\prime}$. Let $A_{1}, A_{2}, \ldots, A_{p + 1}$ be $1$-dimensional subspaces of $W$ such that $A_{i} \subseteq \cap_{V\in\mathcal{V}_{i}} V$ for each $1 \leq i \leq p+1$.

\noindent {\bf Claim 1.} $A_{1}, A_{2}, \ldots, A_{p + 1}$ are distinct.

To the contrary, assume that $A_{i} = A_{j}$ for some $i \neq j$. Then $\mathcal{V}^{\prime} = \mathcal{V}_{i} \cup \mathcal{V}_{j}$ implies that $A_{i} = A_{j} \subseteq \cap_{V\in\mathcal{V^{\prime}}} V$, contradicting to $\dim(\cap_{V\in\mathcal{V^{\prime}}} V) = 0$.

\noindent {\bf Claim 2.} Every $V_{i}$ contains exactly $p$ of the $A_{i}$'s.

Clearly, every $V_{i}$ contains at least $p$ of the $A_{j}$'s since $A_{j} \subseteq V_{i}$ for $i \neq j$. Suppose that $V_{i}$ contains all of the $A_{j}$'s. Then $A_{i} \subseteq V_{i}$ and $A_{i} \subseteq V_{i}\cap (\cap_{V\in\mathcal{V}_{i}} V) = \cap_{V\in\mathcal{V^{\prime}}} V$, a contradiction.

\noindent {\bf Claim 3.} $\dim(V_{i}) = p$ for every $i = 1, 2, \ldots, p + 1.$

For each $1 \leq i \leq p+1$, it is clear that $A_{j} \subseteq V_{i}$ for every $j \neq i$ by Claim $2$. The subspace spanned by $A_{1}$ and $A_{2}$ is $2$-dimensional by Claim $1$.  The subspace spanned by $A_{1}$, $A_{2}$ and $A_{3}$ is $3$-dimensional since $A_{3} \nsubseteq V_{3}$ and $A_{1}\subseteq V_{3}$, $A_{2}\subseteq V_{3}$. Similarly, we conclude that the subspace spanned by $A_{1}$, $A_{2}$, \ldots, $A_{p + 1}$ is $p + 1$-dimensional. Thus, Claim $3$ holds.

Now, the lemma follows from $k \ge \dim(V_{i}) = p$ for every $i = 1, 2, \ldots, p + 1.$
\hfill$\Box$

\vspace{3mm}

\noindent {\bf Lemma 4.2.}
Let $l_{1}$ be a positive integer and $\mathbb{F}_{q}$ be a finite field of order $q$. Let $\mathcal{G}$ be a collection of subspaces of an $n$-dimensional vector space over $\mathbb{F}_{q}$. Suppose that $\dim(\cap_{G \in \mathcal{G}}G) = 0$. Let $V$ be a subspace of the $n$-dimensional vector space with $V \notin \mathcal{G}$ such that $\dim(V \cap G) \ge l_{1}$ for each $G \in\mathcal{G}$. Let $P$ be the subspace spanned by all $G$ in $\mathcal{G}$. Then,
\[
\dim(P \cap V) \ge l_{1} + 1.
\]

\vspace{3mm}

{\bf Proof} $\dim(P \cap V) \ge l_{1}$ since $\dim(V \cap G) \ge l_{1}$. Suppose that $\dim(P \cap V) = l_{1}$. Let $U = P \cap V$. Then
\[
U = P \cap V = \text{span}\{G \in \mathcal{G}\} \cap V \supseteq \text{span}\{G \cap V | G \in \mathcal{G}\}.
\]

Thus, $G \cap V \subseteq U$ for each $G \in \mathcal{G}$. $G \cap V = U$ since $\dim(U) = l_{1}$ and $\dim(G \cap V) = l_{1}$ for each $G \in \mathcal{G}$. Hence $U \subseteq \cap_{G \in \mathcal{G}}G$, a contradiction.
\hfill$\Box$

\vspace{3mm}

\noindent {\bf Lemma 4.3.}
Let $\mathcal{H}$ be a collection of subspaces of an $n$-dimensional vector space over $\mathbb{F}_{q}$. Suppose that $t = |\mathcal{H}| \ge 2$ and $\mathcal{H}$ is a $k$-uniform, intersecting family. Then
\[
\dim(\text{span}\{H | H \in \mathcal{H}\}) \le k + (t - 1)(k - 1).
\]

\vspace{3mm}

{\bf Proof.} We use induction on $t$. It is trivially true for $t = 2$.

Let $t \ge 3$. Suppose that the lemma is true for $t - 1$. Let $\mathcal{H}$ be an arbitrary $k$-uniform intersecting family with $|\mathcal{H}| = t$. Let $\mathcal{G} \subseteq \mathcal{H}$ and $|\mathcal{G}| = t -1$. Clearly $\mathcal{G}$ is a $k$-uniform intersecting family, then
\[
\dim(\text{span}\{G | G \in \mathcal{G}\}) \le k + (t - 2)(k - 1).
\]

Let $\{S\} = \mathcal{H}\backslash \mathcal{G}$. Then,

\[\dim(\text{span}\{H | H \in \mathcal{H}\}) \le \dim(\text{span}\{G | G \in \mathcal{G}\}) + \dim(S) - \dim(\text{span}\{G | G \in \mathcal{G}\}\cap S) \]
\[\le k + (t - 2)(k - 1) + k -1 = k + (t - 1)(k - 1).\]
\hfill$\Box$

\vspace{3mm}

\noindent {\bf Corollary 4.4.}
Let $\mathcal{H}$ be an intersecting family of subspaces of an $n$-dimensional vector space over $\mathbb{F}_{q}$. Suppose that $t = |\mathcal{H}| \ge 2$ and the dimension of every member of $\mathcal{H}$ is at most $k$. Then
\[
\dim(\text{span}\{H | H \in \mathcal{H}\}) \le k + (t - 1)(k - 1).
\]

\vspace{3mm}

We are now ready to prove Theorem 1.15 and the proof is along the same line as the proof for Theorem 1.7.

\vspace{3mm}

{\bf Proof of Theorem 1.15.} Denote $k = \max\{\dim(V_{j}) : 1 \le j \le m\}$. By Theorem 1.13, we may assume that $l_{1} \ge 1$. Also, if there exists $V_{i} \in \mathcal{V}$ such that $\dim(V_{i}) = l_{1}$, then $V_{i} \subseteq \cap_{V\in\mathcal{V}} V$ and the result follows from Theorem 1.13 easily. Thus, we assume that $\dim(V_{i}) \geq l_{1} + 1$ for all $V_{j} \in \mathcal{V}$.

We consider the following cases:

\noindent {\bf Case 1.} $\dim(\cap_{V\in\mathcal{V}} V) = 0$. By Lemma 4.1, there exists a subfamily $\mathcal{V}^{\prime} \subseteq \mathcal{V}$ with $|\mathcal{V}^{\prime}| = k + 1$ such that $\dim(\cap_{V\in\mathcal{V^{\prime}}} V) = 0$. Let
\[
M = \text{span}\{V | V\in\mathcal{V^{\prime}}\}.
\]

Then $\dim(M) \le k + k(k - 1) = k^{2}$ by Corollary 4.4. Since $\dim(V) \ge l_{1} + 1$ for all $V \in \mathcal{V}$, it follows from Lemma 4.2 that
\[ \hspace{32mm} \dim(M \cap V) \ge l_{1} + 1 \quad \text{for each} \quad V \in \mathcal{V}. \hspace{32mm} (4.1)\]
Let $T$ be a given subspace of $M$ such that $\dim(T) = l_{1} + 1$. Define
\[
\mathcal{V}(T) = \{V \in \mathcal{V} : T \subseteq M \cap V\}.
\]
Set $\mathcal{L}^{\prime} = \{l_{2}, l_{3}, \ldots, l_{s}\}$. Then $|\mathcal{L}^{\prime}| = s - 1$. Since $\mathcal{V}$ is $\mathcal{L}$-intersecting family and $\dim(E \cap F)\ge \dim(T) \ge l_{1} + 1$ for any $E, F \in \mathcal{V}(T)$, $\mathcal{V}(T)$ is $\mathcal{L}^{\prime}$-intersecting family. By (4.1), it is easy to see that
\[\hspace{45mm} \mathcal{V} = \cup_{T \subseteq M, \dim(T) = l_{1} + 1}\mathcal{V}(T).  \hspace{45mm} (4.2)\]
Note that for each $T \subseteq M$ with $\dim(T) = l_{1} + 1$, the family
\[
\mathcal{G}(T) = \{V/T : V \in \mathcal{V}(T)\}
\]
is an $\mathcal{L}^{*}$-intersecting family of $W/T$, where $V/T$ is the factor space of $V$ by $T$ and
$\mathcal{L}^{*} = \{l_{2} - l_{1} - 1, l_{3} - l_{1} - 1, \ldots, l_{s} - l_{1} - 1\}$. Clearly, $|\mathcal{G}(T)| = |\mathcal{V}(T)|$. By Theorem 1.13, we have that for each $T \subseteq M$ with $\dim(T) = l_{1} + 1$,
\[
|\mathcal{V}(T)| \le \sum_{j = 0}^{s - 1} \qbinom{n - l_{1} - 1}{j}.
\]

Since $n \ge \log_{q}((q^{s} - 1)\qbinom{k^{2}}{l_{1} + 1} + 1) + l_{1}$, it follows from (4.2) that
\[|\mathcal{V}| \le \sum_{T \subseteq M, \dim(T) = l_{1} + 1}\mathcal{V}(T) \le \qbinom{k^{2}}{l_{1} + 1}\sum_{j = 0}^{s - 1} \qbinom{n - l_{1} - 1}{j}\]
\[\le \qbinom{k^{2}}{l_{1} + 1}\frac{q^{s} - 1}{q^{n - l_{1}} - 1}\sum_{j = 0}^{s - 1} \qbinom{n - l_{1}}{j + 1} \le \sum_{j = 0}^{s} \qbinom{n - l_{1}}{j}.\]

\noindent {\bf Case 2.} $\dim(\cap_{V\in\mathcal{V}} V) \not= 0$. For this case, if $\dim(\cap_{V \in \mathcal{V}} V) \ge l_{1}$, then the result follows from easily from Theorem 1.13. Assume that $0 < \dim(\cap_{V\in\mathcal{V}} V) = t < l_{1}$ and let $T = \cap_{V\in\mathcal{V}} V$. Then
\[
\mathcal{G} = \{V/T : V \in \mathcal{V}\}
\]
is an $\mathcal{L^{\prime}}$-intersecting family of $W/T$, where $\mathcal{L^{\prime}} = \{l_{1} - t, l_{2} - t, \ldots, l_{s} - t\}$ with $l_{1} - t > 0$. By Case 1, we obtain
\[
|\mathcal{V}| = |\mathcal{G}| \le \sum_{j = 0}^{s}\qbinom{(n - t) - (l_{1} - t)}{j} =  \sum_{j = 0}^{s}\qbinom{n - l_{1}}{j}.
\]
\hfill$\Box$

\vspace{3mm}

Similarly, one can prove Theorem 1.16 by using Theorem 1.14 instead of Theorem 1.13.

\end{document}